\documentstyle{amsppt}
\tolerance 3000
\pagewidth{5.5in}
\vsize7.0in
\magnification=\magstep1
\widestnumber \key{AAAAAAA}
\topmatter
\author Alex Iosevich    
\endauthor 
\thanks Research at MSRI supported in part by NSF grant DMS97-06825
\endthanks 
\address Department of Mathematics Georgetown University Washington D.C. 20057 \newline 
email: iosevich $\@$math.georgetown.edu \endaddress 
\title Maximal averages over flat radial hypersurfaces   
\endtitle 
\endtopmatter 

Let $A_tf(x)=\int_S f(x-ty) d\sigma(y)$, where $S$ is a smooth compact
hypersurface in ${\Bbb R}^n$ and $d\sigma$ denotes the
Lebesgue measure on $S$. Let ${\Cal A}f(x)=\sup_{t>0}|A_tf(x)|$. 
If the hypersurface $S$ has non-vanishing Gaussian curvature, then 
$$ {||{\Cal A}f||}_{L^p({\Bbb R}^n)} \leq C_p {||f||}_{L^p({\Bbb R}^n)}, \ \ 
f \in {\Cal S}({\Bbb R}^n), \tag*$$ for $p>\frac{n}{n-1}$. Moreover, the 
result is sharp. See \cite{St76}, \cite{Gr82}. 

If the hypersurface $S$ is convex and the order of contact with every 
tangent line is finite, the optimal exponents for the inequality $(*)$ 
are known in ${\Bbb R}^3$, (see \cite{IoSaSe97}), and in any dimension 
in the range $p>2$, (see \cite{IoSa96}). More precisely, the result in 
the range $p>2$ is the following. 
\proclaim{Theorem 1 (\cite{IoSa96})} Let $S$ be a smooth convex compact 
finite type hypersurface, in the sense that the order of contact with 
every tangent line is finite. Then for $p>2$, the following condition
is necessary and sufficient for the maximal inequality $(*)$. 
$$ {(d(x, {\Cal H}))}^{-1} \in L^{\frac{1}{p}}(S), \tag1$$ for every 
tangent hyperplane ${\Cal H}$ not passing through the origin, 
where $d(x, {\Cal H})$ denotes the 
distance from a point $x \in S$ to the tangent hyperplane ${\Cal H}$. 
\endproclaim 

In fact, the condition $(1)$ is a necessary condition for any smooth 
compact hypersurface in ${\Bbb R}^n$. See \cite{IoSa96}, Theorem 2. 

In this paper we shall consider convex radial hypersurfaces of the form 
$$S=\{x \in B: x_n=\gamma(|x'|)+1\},\tag2$$ where $B$ is a ball centered at the 
origin, $x=(x',x_n)$, $\gamma$ is convex, $\gamma$ , $\gamma''$ increasing,  
$\gamma(0)=\gamma'(0)=0$, and 
$\gamma$ is allowed to vanish of infinite order. 

If $\gamma''$ does vanish of infinite order, the condition $(1)$ cannot 
hold for any $p<\infty$. Since the condition $(1)$ is necessary by Theorem 1 above, our only hope is to look for an inequality 
of the form 
$$ {||{\Cal A}f||}_{L^{\Phi}({\Bbb R}^n)} \leq C_{\Phi} {||f||}_{L^{\Phi}
({\Bbb R}^n)}, \tag3$$ where $L^{\Phi}({\Bbb R}^n)$ is an Orlicz space, 
near $L^{\infty}({\Bbb R}^n)$, associated to a Young function $\Phi$, with
the norm given by 
$$ {||f||}_{\Phi}=\inf \left\{s>0: \int \Phi\left( \frac{|f(x)|}{s} \right)dx 
\leq 1\right\}. \tag4$$ 

The following result was proved in \cite{Bak95}. 
\proclaim{Theorem 2} Let $S$ be as in $(2)$ with $n=3$. Assume that for 
each $\lambda>1$ 
$$ \frac{\gamma'(\lambda t)}{\gamma'(t)} \ \ \text{is non-decreasing for} \ \ 
t>0. \tag5$$ 

Put $G(t)=t^2\gamma'(t)$. For $\beta>1$ and $d>0$ let $\phi:[0, \infty) 
\rightarrow [0, \infty)$ be a non-decreasing function such that 
$\phi(t)=t^{-1}{[G(t^{-d})]}^{-\beta}$ if $t$ is sufficiently large, 
$\phi(t)>0$ if $t>1$, and $\phi(t)=0$ if $0 \leq t \leq 1$. Let 
$\Phi(u)=\int_{0}^{u} \phi(t)dt$. Then for every $d>\frac{1}{2}$ there 
exists a constant $C$ such that the estimate $(3)$ holds. 
\endproclaim 

The examples show (see \cite{Bak95}, Example 3.3) that Theorem 2 is sharp 
for some surfaces, for example if $\gamma(s)=e^{-\frac{1}{s^b}}$, $b>0$, but 
not for others, for example if $\gamma(s)=s^m$. 

In this paper we shall give a set of simple sufficient conditions for the 
inequality $(3)$ for some classes of Orlicz functions $\Phi$. We will 
show that our result is sharp for a wide class of both finite type and 
infinite type $\gamma$'s. 

{\bf Acknowledgements:} The author wishes to thank Jim Wright for teaching him the technique needed to prove the three dimensional case of Lemma 5 below. 

\head Assumptions on $\Phi$ \endhead
\vskip.125in 

Assume that $\Phi$ is a Young function such that $\Phi(s)=\int_{0}^{s} 
\phi(t)dt$, where $\phi:[0, \infty) \rightarrow [0, \infty)$ is a 
non-decreasing function such that $\phi(t)=0$ for $0 \leq t \leq 1$, and 
$\phi(t)>0$ for $t>1$. Assume that there exist constants $c>1$, $C_0$, and 
$C_1$ such that 
$$ \int_{1}^{u} \frac{\phi(t)}{t^r}dt \leq C_0 \frac{\phi(u)}{u^{r-1}} \ \ 
\text{for} \ \ u>1, \tag6$$ and for every $\lambda>1$ 
$$ C_1 \frac{\phi(\lambda t)}{\phi(t)} \ge \phi(\lambda) \ \ \text{for} \ \ 
t \ge c. \tag7$$ 

Our main reason for making these assumptions about $\Phi$ is the following
generalization of the Marcienkiewicz interpolation theorem due to Bak. See 
\cite{Bak95}, Lemma 1.1. 
\proclaim{Lemma 3} Let $r \in [1,\infty)$. Suppose that the operator $T$ 
is simultaneously weak type $(1,1)$ and $(\infty, \infty)$, namely there
exist constants $A, B>0$ such that 
$$ \mu( \{x: |Tf(x)|>t \}) \leq {\left(\frac{A{||f||}_r}{t} \right)}^r \ \ 
\forall \ t>0, \tag8$$ 
$$ {||Tf||}_{\infty} \leq B{||f||}_{\infty}. \tag9$$ 

Suppose that $\Phi$ satisfies the assumptions above. Then there exists a
constant $C=C(\Phi, r)$ depending only on $\Phi$ and $r$ such that 
$$ {||Tf||}_{\Phi} \leq CB\Phi^{-1}({(A/B)}^r) {||f||}_{\Phi}. \tag10$$ 
\endproclaim 

\remark{Remark} Lemma 3 has the following interesting  consequence. Let 
${\Cal A}f(x)=\sup_{t>0} \int f(x-t(s,s^m+1)) \psi(s)ds$, $m>2$, where $\psi$ 
is a smooth cutoff function, and let ${\Cal A}^kf(x)$ denote the same 
operator with $s$ localized to the interval $[2^{-k}, 2^{-k+1}]$. It 
was proved in \cite{I94} that ${\Cal A}^k: L^p({\Bbb R}^2) \rightarrow
L^p({\Bbb R}^2)$, $p>2$, with norm $C2^{-k} 2^{\frac{mk}{p}}$. Let 
$\Phi_{p, \alpha}(t)=t^p \log^{\alpha}(t)$. It follows by Lemma 3 that 
${\Cal A}: L^{\Phi_{p, \alpha}}({\Bbb R}^2) \rightarrow 
L^{\Phi_{p, \alpha}}({\Bbb R}^2)$ if 
$p=m$ and $\alpha>m$. \endremark 

\head Statement of results \endhead 
\vskip.125in 

Our main results are the following. 
\proclaim{Theorem 4} Let $S$ be as in $(2)$. Let $n \ge 3$. Suppose that $\Phi$ 
satisfies the conditions $(6)$ and $(7)$ above. Suppose that $\lim_{t  
\rightarrow 0} \Phi(t)/t^2=0$. Then the estimate $(3)$ holds if  
$$\sum_{j=0}^{\infty} 2^{-j(n-1)} \Phi^{-1} \left( \frac{1}{\gamma(2^{-j})} 
\right)<\infty. \tag11$$.
\endproclaim 

The main technical result involved in the proof of Theorem 4 is the following
version of the standard stationary phase estimates. 
\proclaim{Lemma 5} Let $n \ge 3$. Let 
$$ F_j(\xi)=\int_{\{y: 1 \leq |y| \leq 2\}} e^{i(\langle y, \xi' \rangle +
\xi_n \gamma_j(|y|))} e^{i\frac{\xi_n}{\gamma(2^{-j})}} dy, \tag12$$ 
with $\gamma_j(s)=\frac{\gamma(2^{-j}s)}
{\gamma(2^{-j})}$, where $\gamma$ is as in $(2)$. Then 
$$ |F_j(\xi)| \leq C{(1+|\xi|)}^{-1}, \tag13$$   
where $C$ is independent of $j$ and $\gamma$.  

Moreover, if $|F_j(\xi)|$ is replaced by $|\nabla F_j(\xi)|$ then the 
estimate $(13)$ still holds with $C$ on the right-hand side
replaced by $\frac{C}{\gamma(2^{-j})}$. 
\endproclaim 

The main technical result used in the proof of Theorem 2 is the following. 
See \cite{Bak95}, Theorem 2.1. 
\proclaim{Lemma 6} Let $\chi \in C^1_{0}([0, \infty)])$ be a non-negative
function that is compactly supported in the interval $(a, \infty)$, 
where $a>0$. Let $n=3$ and let $S$ be as in $(2)$ where $\gamma$ satisfies
the condition of Theorem 2. Let $F_S(\chi)(\xi)$ denote $F_0(\xi)$ in 
Lemma 5 with $\chi(|y|)$ in place of the characteristic function of 
the annulus $\{y: 1 \leq |y| \leq 2\}$. 

Then for every multi-index  $\alpha$ with $|\alpha| \leq 1$ there exists a 
constant $C$ independent of $a$, $\xi$, and $\chi$ such that 
$$ |{(\partial / \partial \xi)}^{\alpha} F_S(\chi)(\xi)| \leq 
C C_{\chi} \frac{a}{\sqrt{\gamma'(a) \gamma'(a/2)}} {(1+|\xi|)}^{-1}, 
\tag14$$ where $C_{\chi} \leq {||\chi||}_{\infty}+{||\chi'||}_1$ if $\alpha=0$,
and $C_{\chi} \leq {||\chi||}_{\infty}+{||\chi||}_1+{||\chi'||}_1$ if 
$\alpha=1$. 
\endproclaim  

\head Main idea \endhead 
\vskip.125in 

The point is that even though a higher dimensional analog of Lemma 6 may be
difficult to obtain, we get around the problem by using Lemma 5. We have
to settle for the uniform decay of order 
$\max \{ -\frac{n-2}{2}, -1 \}$ instead of
$-\frac{n-1}{2}$, but this is enough in dimension $n \ge 4$ as we shall see 
below. The idea is, roughly speaking, the following. We are trying to prove 
$L^{\Phi} \rightarrow L^{\Phi}$ estimates for maximal operators associated to 
radial convex surfaces. If the surface is infinitely flat, then Theorem 2 
in \cite{IoSa96} implies that 
$L^p \rightarrow L^p$ estimates are not possible for $p<\infty$. 
So we are looking for $L^{\Phi}\rightarrow L^{\Phi}$ estimates where $L^{\Phi}$
is very close to $L^{\infty}$, so interpolating between $L^2$ and $L^{\infty}$
in the right way should do the trick. However, in order to obtain $L^2$ 
boundedness of the maximal operator, we only need decay $-\frac{1}{2}-
\epsilon$, $\epsilon>0$. If $n \ge 4$, then $\frac{n-2}{2}>\frac{1}{2}$, so we
should be alright. If $n=3$ a bit more integration by parts will be required.   

\head Plan \endhead 
\vskip.125in 

The rest of the paper is organized as follows. In the next section we shall prove 
Theorem 4 assuming Lemma 5. In the following section we shall prove Lemma 5.
In the final section of the paper we shall discuss the sharpness of Theorem 4
and give some examples.  

\head Proof of Theorem 4 \endhead 
\vskip.125in 

Let $A^j_tf(x)=\int f(x'-ty, x_n-t(\gamma(|y|)+1)) \psi_0(y)dy$, where
$\psi_0$ is a smooth cutoff function supported in $[1,2]$, such that 
$\sum_j \psi(2^js) \equiv 1$. Let $\tau_jf(x)=f(2^{-j}x', \gamma(2^{-j})x_n)$.
Making a change of variables we see that 
$$ A^j_tf(x)=2^{-j(n-1)} \tau_j^{-1} B^j_t \tau_jf(x), \tag15$$ where 
$$ B^j_tf(x)=\int f(x'-ty, x_n-t(\gamma_j(|y|)+1/\gamma(2^{-j})) \psi_0(y)dy. 
\tag16$$  

We shall prove that 
$$ \sup_{t>0}B^j_t: L^2({\Bbb R}^n) \rightarrow L^2({\Bbb R}^n) \ \ \text{ with
norm} \ \  {\left( \frac{1}{\gamma(2^{-j})} \right)}^{\frac{1}{2}}.\tag17$$ 

By interpolating with the trivial estimate ${||\sup_{t>0}B^j_tf||}_{\infty} 
\leq C{||f||}_{\infty}$ using Lemma 3, we shall conclude that 
$$ \sup_{t>0}B^j_t: L^{\Phi}({\Bbb R}^n) \rightarrow L^{\Phi}({\Bbb R}^n) \ \ 
\text{with norm} \ \ \Phi^{-1}\left(\frac{1}{\gamma(2^{-j})} \right). \tag18$$

Since the $L^p$ norms of $\tau_j$ and $\tau_j^{-1}$ are reciprocals of 
each other, it follows that ${\Cal A}: L^{\Phi}({\Bbb R}^n) \rightarrow 
L^{\Phi}({\Bbb R}^n)$ if 
$$ \sum_{j=0}^{\infty} 2^{-j(n-1)} \Phi^{-1}\left( \frac{1}{\gamma(2^{-j})} 
\right)<\infty. \tag19$$ 

So it remains to prove $(18)$. The proof follows from the standard 
Sobolev imbedding theorem type argument. See for example \cite{St76}. 
We shall use the following version which follows from the proof of Theorem 
15 in \cite{IoSa96}. See also, for example, \cite{CoMa86}, \cite{MaRi95}. 
\proclaim{Lemma 7} Suppose that $\tau$ is the Lebesgue measure on the 
hypersurface $S$ supported in an ellipsoid with eccentricities  
$(1, \dots, 1, R)$. Suppose that  $|\hat{\tau}(\xi)| \leq C$ and 
$\max \{|x|: x \in supp(\tau) \} \leq 10R$. Suppose that 
$$ { \left( \int_{1}^{2} { \left| \hat{\tau}(t\xi) \right|}^2 dt \right)}^
{ \frac{1}{2}} \leq C {(1+|\xi|)}^{-\frac{1}{2}-\epsilon}, and \tag20$$
$$ {\left( \int_{1}^{2} {\left| \nabla \hat{\tau}(t\xi) \right|}^2 dt 
\right)}^{\frac{1}{2}} \leq CR {(1+|\xi|)}^{-\frac{1}{2}-\epsilon} \tag21$$
for some $\epsilon>0$. Let $\hat{\tau_t}(\xi)=\hat{\tau}(t\xi)$. Let
${\Cal M}f(x)=\sup_{t>0}|f*\tau_t(x)|$. Then 
$$ {||{\Cal M}f||}_2 \leq 100C \sqrt{R}{||f||}_2. \tag22$$ 
\endproclaim 

Application of Lemma 7 immediately yields $(17)$ since by Lemma 5 $C$  
is a universal constant and $R \leq \frac{C}{\gamma(2^{-j})}$. This completes
the proof of Theorem 4. 

\head Proof of Lemma 5 \endhead 
\vskip.125in 

We must show that  
$$ |F_j(\xi)|=\left|\int_{\{y: 1 \leq |y| \leq 2\}}  
e^{i(\langle y, \xi' \rangle+\xi_n \gamma_j(|y|))} 
e^{i\frac{\xi_n}{\gamma(2^{-j})}} dy\right| 
\leq C{|\xi|}^{-1}, \tag23$$ with $C$ independent of $\gamma$ and $j$. 

Our plan is as follows. We will first show that if either 
$|\xi'| \approx |\xi_n|$, or $|\xi'| >> |\xi_n|$, then 
$|F_j(\xi)| \leq C{(1+|\xi|)}^{-\frac{n-2}{2}}$. If $|\xi_n| >> |\xi'|$, 
we will show that $|F_j(\xi)| \leq C{(1+|\xi_n|)}^{-1}$. This will 
complete the proof since $\frac{n-2}{2} \ge 1$ if $n \ge 4$.

Going into polar coordinates and applying stationary phase, we get  
$$ e^{i \frac{\xi_n}{\gamma(2^{-j})}} 
\int_{1}^{2} e^{i\xi_n \gamma_j(r)} r^{n-2} dr \int_{S^{n-2}} 
e^{ir \langle \xi', \omega \rangle} d\omega. \tag24$$ 

Since the Gaussian curvature on $S^{n-2}$ does not vanish, it is a 
classical result that 
$$ \left| \int_{S^{n-2}} e^{i \langle \xi', \omega \rangle} d\omega 
\right| \leq C{(1+|\xi'|)}^{-\frac{n-2}{2}}. \tag25$$ 

It follows that $|F_j(\xi)| \leq 
C{(1+|\xi|)}^{-\frac{n-2}{2}}$ if either $|\xi'| >> |\xi_n|$ or 
$|\xi'| \approx |\xi_n|$. If $|\xi_n| >> |\xi'|$, let $h(r)=\xi_n \gamma_j(r)
-r\langle \xi', \omega \rangle$. Since $\gamma$ is convex, it follows
that $|h'(r)| \ge |\xi_n|-|\xi'|$. Since $|\xi_n| >> |\xi'|$, it follows
by the Van Der Corput Lemma that the expression in $(24)$ is bounded by 
$C \frac{1}{|\xi|}$. 

The estimate for $\nabla F_j$ follows in the same way if we observe that 
the derivative with respect to $\xi_n$ brings down a factor of 
$\gamma_j(r)+\frac{1}{\gamma(2^{-j})}$,
and $\gamma_j(r) + \frac{1}{\gamma(2^{-j})} \leq 2 \frac{1}{\gamma(2^{-j})}$. 
This completes the proof of Lemma 5 if $n \ge 4$. 

To prove the three dimensional case we go into polar coordinates, integrate in the angular variables and use the well known asymptotics for the Fourier transform of the Lebesgue measure on the circle to obtain 
$$ \int e^{i\phi(r)} r b(rA) \psi_0(r)dr, \tag26$$ where $A=|\xi'|$, $\lambda=\xi_n$, $b$ is a symbol of order $-\frac{1}{2}$, $\psi_0$ is as above, and $\phi(r)=rA-\gamma_j(r)\lambda$. 

Let 
$$ G(r)=\int_{r}^{2} e^{i\phi(s)} ds, \tag27$$ so the integral in $(26)$ becomes  
$$ \int G'(r) r b(rA)\psi_0(r)dr. \tag28$$ 

Integrating by parts we get 
$$ \int G(r) {(r b(rA) \psi_0(r))}' dr. \tag29$$ 

Let $r_0$ be defined by the relation $\gamma_j'(r_0)=\frac{A}{2\lambda}$.
We have $|\phi''(s)| \ge |\gamma_j''(s) \lambda| \ge 
|\gamma_j'(s) \lambda| \ge |\gamma_j'(r) \lambda|$. If $r_0<r$ this quantity 
is bounded below by $C|A|$ and the Van der Corput lemma gives the decay  
$C{|A|}^{-\frac{1}{2}}$ for $G(r)$. Using the 
fact $b$ is a symbol of order $-\frac{1}{2}$ we see that $(29)$ is 
bounded by $C{|A|}^{-1}$, $|A|$ large. This handles the 
case $|\lambda| \leq C|A|$ and $r \leq r_0$. 

On the other hand, $|\phi'(s)|=|A-\gamma_j'(s)\lambda|$. Split up the integral 
that defines $G(r)$ into two pieces: $s \in [r,r_0]$ and $s \in [r_0,2]$. The
second integral was just handled above. In the first integral 
$|\phi'(s)| \ge |\phi'(r_0)| \ge C|A|$. The Van der Corput Lemma yields decay 
$\frac{C}{|A|}$. Taking the properties of the symbol $b$ into account, 
as before, we get the decay $\frac{C}{|A|} {|A|}^{-\frac{1}{2}}$. 
This takes care of the case $|\lambda| \leq C|A|$ and $r \ge r_0$.  

If $|\lambda|>>|A|$, $|\phi'(s)| \ge C|\lambda|$ and the Van der Corput 
lemma yields the decay $\frac{C}{|\lambda|}$ for $(29)$. This completes
the proof of the three dimensional case. 

\head Examples \endhead 
\vskip.125in 

\proclaim{Example 1} Let $\gamma(s)=s^m$, $m \ge 2(n-1)$, and $\Phi(t)=t^p$. 
Theorem 4 yields boundedness for $p>\frac{m}{n-1}$. 
This is sharp by Theorem 1. \endproclaim 

\proclaim{Example 2} Let $\gamma(s)=s^m$, $m \ge 2(n-1)$, and 
$\Phi_{p,\alpha}(s)=s^p \log^{\alpha}(s)$. Then Theorem 4 yields 
boundedness for $p=\frac{m}{n-1}$ and $\alpha>\frac{m}{n-1}$. 
\endproclaim  

\proclaim{Example 3} Let $\gamma(s)=e^{-\frac{1}{s^{\alpha}}}$, $\alpha>0$, and 
$\Phi(t)=e^{t^{\beta}}$, $\beta>0$. Then Theorem 4 tells us that the maximal 
operator is bounded if $\alpha<\beta(n-1)$. Testing $A_tf(x)$ against 
$h_p(x)=\Phi^{-1}\left(\frac{1}{|x_n|} \right) 
\frac{1}{\log \left( \frac{1}{|x_n|}
\right)} \chi_B(x)$, where $\chi_B$ is the characteristic function of the ball
of radius $\frac{1}{2}$ centered at the origin, shows that this result 
is sharp. The same procedure establishes sharpness of the estimate given 
in Example 2. \endproclaim 

In fact, testing $A_tf(x)$ against $h_p(x)$ shows that the summation 
condition of Theorem 4 is pretty close to being sharp. It is not hard to see 
that, at least up to a $\log$ factor, ${\Cal A}$ bounded on 
$L^{\Phi}({\Bbb R}^n)$ implies that 
$$\int_{\{y: |y| \leq 2\}} 
\Phi^{-1}\left( \frac{1}{\gamma(|y|)} \right) dy<\infty.\tag30$$ This would 
literally follow, without the log factor, from the proof of Theorem 2 of 
\cite{IoSa96} 
if we assumed, in addition, that $\Phi(ab) \ge \Phi(a)\Phi(b)$, for every
$a,b>0$. 

The condition $(30)$ is equivalent (after making a change of variables 
and going into polar coordinates) to 
$$ \sum_{j=0}^{\infty} 2^{-j(n-1)} \int_{1}^{2} \Phi^{-1} \left(
\frac{1}{\gamma(2^{-j}r)} \right) r^{n-2}dr<\infty. \tag31$$ 

The expression $(31)$ is equivalent to the summation condition of Theorem 4
if $\gamma$ does not vanish to infinite order. If $\gamma$ vanishes to 
infinite order, the two condition are still often equivalent, as in 
the Example 3 above. 

\remark{Remark} It would be interesting to extend the results of this paper to a more general class of hypersurfaces. For example, one could consider hypersurfaces of the form $S=\{x \in {\Bbb R}^n: x_n=\gamma(\phi(x'))+1\}$ where $\gamma$ is as above and $\phi$ is a smooth convex finite type function. Some recent results (see e.g. \cite{IoSa97}, \cite{IoSaSe98}, and \cite{WWZ97}) suggest that such an analysis should be possible. We shall address this issue in a subsequent paper (\cite{I98}. More generally, a bigger challange would be to consider a hypersurface of the form $S=\{x \in {\Bbb R}^n: x_n=G(x')+1\}$, where $G$ is a smooth function of $n-1$ variables that vanishes of infinite order at the origin. At the moment, obtaining sharp Orlicz estimates, even in the case where the determinant of the Hessian matrix of $G$ only vanishes at the origin, does not seem accessible. 
\endremark 

\newpage 

\head References \endhead 
\vskip.25in 

\ref \key Bak95 \by J.-G. Bak \paper Averages over surfaces with infinitely 
flat points \jour J. Func. Anal. \yr 1995 \vol 129 \endref 

\ref \key B86 \by J. Bourgain \paper Averages in the plane over convex curves
and maximal operators \jour J. Analyse Math. \vol 47 \yr 1986 \pages 69-85
\endref 

\ref \key CoMa86 \by M. Cowling and G. Mauceri \paper Inequalities for some
maximal functions II \jour Trans. of the A.M.S. \yr 1986 \vol 296 \endref 

\ref \key Gr82 \by A. Greenleaf \pages 519-537 
\paper Principal curvature and harmonic analysis 
\yr 1982 \vol30 \jour Indiana Math J. \endref 

\ref \key I94 \by A. Iosevich 
\paper Maximal operators associated to families of flat curves in the plane
\yr 1994 \jour Duke Math J. \vol 76 \endref 

\ref \key I98 \by A. Iosevich \paper Averages over convex infinite type hypersurfaces \jour (in preparation) \endref 

\ref \key IoSa96 \by A. Iosevich and E. Sawyer 
\paper Oscillatory integrals and maximal averaging operators associated to 
homogeneous hypersurfaces \yr 1996 \jour Duke Math J. \vol 82 \endref 

\ref \key IoSa97 \by A. Iosevich and E. Sawyer 
\paper Maximal averages over surfaces  
\jour Adv. in Math. \vol 132 \yr 1997 \endref 

\ref \key IoSaSe98 \by A. Iosevich, E. Sawyer, and A. Seeger \paper Averages
over convex surfaces of finite type in ${\Bbb R}^3$ \jour (preprint)
\yr 1998 \endref 

\ref \key MaRi95 \by G. Marletta and F. Ricci \paper Two parameter 
maximal functions associated with homogeneous surfaces in ${\Bbb R}^n$ 
\jour preprint \yr 1995 \endref 

\ref \key St76 \by E. Stein \paper Maximal functions: spherical means 
\yr 1976 \jour Proc. Nat. Acad. Sci. U.S.A. \vol73 \endref 

\ref \key St93 \by E. Stein \paper Harmonic Analysis 
\jour Princeton Univ. Press \yr 1993 \endref

\ref \key WWZ97 \by J. Wright, S. Wainger, and S. Ziesler \paper (preprint) \endref 

\enddocument